\documentclass[10pt]{amsart}
\usepackage{amsmath, amsthm, amscd, amsfonts}

\newtheorem{theorem}{Theorem}[section]

\newtheorem{proposition}[theorem]{Proposition}
\newtheorem{corollary}[theorem]{Corollary}
\theoremstyle{definition}

\theoremstyle{remark}

\numberwithin{equation}{section}
\newfont{\kh}{msbm10}
\newcommand{\R}{\mbox{\kh R}}
\newcommand{\N}{\mbox{\kh N}}

\begin{document}
\title{Stability of Generalized Jensen Equation on Restricted Domains}
\author{S.-M. Jung}
\address{Soon-Mo Jung: Mathematics Section, College of Science and Technology, Hong-Ik
University, 339-701 Chochiwon, Korea}
\email{smjung@wow.hongik.ac.kr}
\author{M. S. Moslehian}
\address{Mohammad Sal Moslehian: Department of
Mathematics, Ferdowsi University, P. O. Box 1159, Mashhad 91775,
Iran}
\email{moslehian@ferdowsi.um.ac.ir}
\author{P. K. Sahoo}
\address{Prasanna K. Sahoo: Department of Mathematics, University of Louisville, Louisville, Kentucky
40292, USA}
\email{sahoo@louisville.edu} 

\begin{abstract}
In this paper, we establish the conditional Hyers-Ulam-Rassias
stability of the generalized  Jensen functional
equation $r \, f \left (\frac{sx+ty}{r} \right ) = s \, g(x) + t\, h(y)$ on various
restricted domains such as inside balls, outside balls, and
punctured spaces.
In addition, we prove the orthogonal stability
of this equation and study orthogonally generalized Jensen
mappings on Balls in inner product spaces.
\end{abstract}
\maketitle

\section{Introduction}
It is known that the problem of stability of functional equations
originated from the following question of
Ulam \cite{ULA} posed in 1940:
Given an approximately linear mapping $f$, when does a linear
mapping $T$ estimating $f$ exist? In the next year,
Hyers \cite{HYE} gave a partial affirmative answer to the
question of Ulam in the context of Banach spaces.
In 1978, Rassias \cite{RAS1} extended the theorem of Hyers
by considering the unbounded Cauchy difference
$\| f(x+y) - f(x) - f(y) \| \leq \varepsilon( \|x\|^p + \|y\|^p)$,
($\varepsilon>0$, $p\in[0,1)$).
The result of Rassias has provided a lot of influence in the
development of what we now call {\it Hyers--Ulam--Rassias
stability} of functional equations. In 1992, a generalization of
Rassias' theorem, the so-called generalized Hyers--Ulam--Rassias
stability, was obtained by G\u avruta \cite{GAV}. During the last
decades several stability problems of functional equations have
been investigated in the spirit of Hyers--Ulam--Rassias-G\u
avruta. The reader is referred to \cite{CZE,H-I-R,JUNG1,MOS1,RAS2}
and references therein for more detailed information on stability
of functional equations.

The generalized Jensen equation is $r \, f(\frac{sx+ty}{r}) = s \, f(x)
+ t \, f(y)$ where $f$ is a mapping between linear spaces and $r, s,
t$ are given constant positive integers. It is easy to see that a
mapping $f : {\mathcal X} \to {\mathcal Y}$ between linear spaces
with $f(0)=0$ satisfies the generalized Jensen equation for all
$x, y \in {\mathcal X}$ if and only if it is additive; cf.
\cite{B-O-P-P} (see also Theorem 6 of \cite{P-V} in the case of
Jensen's equation).

Stability of (generalized) Jensen equation has been studied at
first by Kominek \cite{KOM} and then by several other
mathematicians; cf. \cite{C-R, JUNG, L-J}. In this paper,
using the ``direct method'',  we establish the
Hyers--Ulam--Rassias stability of the generalized Jensen
functional equation of Pexider type and the conditional stability
on some certain restricted domains. Throughout the paper,
${\mathcal X}$ denotes a linear space and ${\mathcal Y}$
represents a Banach space. In addition, we assume $r, s, t$ to be
constant positive integers.

\section{Stability of generalized Jensen equation}

The following theorem is a simple generalization of Theorem 2.1 of
\cite{B-O-P-P}. It can also be regarded as an extension of the
main theorem of \cite{J-S-K}.

\begin{theorem}
Let $f, g, h : {\mathcal X} \to {\mathcal Y}$ be mappings with
$f(0) = g(0) = h(0) = 0$ for which there exists a function
$\varphi : {\mathcal X}\times {\mathcal X} \to [0, \infty)$
satisfying
\begin{eqnarray*}
\widetilde{\varphi}(x, y) : = \frac{1}{2r} \sum_{n=0}^{\infty}
\frac{1}{2^n} \left[ \varphi\!\left( 2^n\frac{r}{s}x,
2^n\frac{r}{t}y \right) +
              \varphi\!\left( 2^n\frac{r}{s}x, 0 \right) +
              \varphi\!\left( 0, 2^n\frac{r}{t}y \right)
       \right]
< \infty ,
\end{eqnarray*}
and
\begin{eqnarray}\label{Jen}
\left\| r \, f\!\left( \frac{sx+ty}{r} \right) - s \, g(x) - t \, h(y) \right\|
\leq \varphi(x,y),
\end{eqnarray}
for all $x, y\in {\mathcal X}$.
Then there exists a unique additive map $T: {\mathcal X}\to {\mathcal Y}$
such that
\begin{eqnarray*}
\| f(x) - T(x) \| \leq \widetilde{\varphi}(x,x)
\end{eqnarray*}
\begin{eqnarray*}
\| g(x) - T(x) \| \leq \frac{1}{s}\, \varphi(x,0) +
\frac{r}{s}\, \widetilde{\varphi}\!\left( \frac{s}{r}x, \frac{s}{r}x \right),
\end{eqnarray*}
\begin{eqnarray*}
\| h(x) - T(x) \| \leq \frac{1}{t}\, \varphi(0,x) +
\frac{r}{t}\, \widetilde{\varphi}\!\left( \frac{t}{r}x, \frac{t}{r}x \right),
\end{eqnarray*}
for all $x\in {\mathcal X}$.
Furthermore, if the mapping $\mu \mapsto f(\mu x)$ is continuous
for each fixed $x\in {\mathcal X}$, then the additive map $T$
is $\R$-linear.
\end{theorem}

\begin{proof}
Setting $y=0$ in (\ref{Jen}) we get
\begin{eqnarray}\label{f-g}
\left\| r\,  f\!\left( \frac{s}{r}x \right) - s \, g(x) \right\|
\leq \varphi(x,0)
\end{eqnarray}
for all $x\in {\mathcal X}$.
Setting $x=0$ in (\ref{Jen}) we obtain
\begin{eqnarray*}
\left\| r \, f\!\left( \frac{t}{r}y \right) - t \, h(y) \right\|
\leq \varphi(0,y)
\end{eqnarray*}
for all $y\in {\mathcal X}$.
Then
\begin{eqnarray*}
\left\| r f\!\left( \frac{sx+ty}{r} \right) \right . &-&
      \left.  r f\!\left( \frac{s}{r}x \right) -
        r f\!\left( \frac{t}{r}y \right) \right\| \\
& \leq & \left\| r f\!\left( \frac{sx+ty}{r} \right) - s \, g(x) - t \, h(y) \right\| \\
&    &+ \,\left\| r f\!\left( \frac{s}{r}x  \right) - s \, g(x) \right\|
         +  \left\| r f\!\left( \frac{t}{r}y \right) - t \, h(y) \right\| \\
& \leq & \varphi(x,y) + \varphi(x,0) + \varphi(0,y)
\end{eqnarray*}
for all $x, y\in {\mathcal X}$.
Replacing $x$ and $y$ by $\frac{r}{s}x$ and $\frac{r}{t}y$, respectively,
we have
\begin{eqnarray*}
\| f(x+y) - f(x) - f(y) \| \leq \frac{1}{r}
\left[ \varphi\!\left( \frac{r}{s}x, \frac{r}{t}y \right) +
       \varphi\!\left( \frac{r}{s}x, 0 \right) +
       \varphi\!\left( 0, \frac{r}{t}y \right)
\right].
\end{eqnarray*}
By the well-known theorem of G\u avruta \cite{GAV} there exists a unique
additive mapping $T: {\mathcal X}\to {\mathcal Y}$ given by
$T(x) := \displaystyle{\lim_{n\to\infty}} 2^{-n}f(2^nx)$ such that
\begin{eqnarray}\label{f-T}
\|f(x) - T(x)\| \leq \widetilde{\varphi}(x, x)
\end{eqnarray}
for all $x \in {\mathcal X}$.
Since $T$ is additive, we have $T(\alpha x) = \alpha T(x)$ for all
rational numbers $\alpha$ and $x \in {\mathcal X}$.
It follows from (\ref{f-g}) and (\ref{f-T}) that
\begin{eqnarray*}
\| g(x) - T(x) \|
& \leq & \left\| g(x) - \frac{r}{s} \, f\!\left( \frac{s}{r} x \right) \right\| +
         \left\| \frac{r}{s} \,  f\!\left( \frac{s}{r} x \right) - T(x) \right\| \\
& \leq & \frac{1}{s} \, \varphi(x,0) +
         \frac{r}{s} \, \widetilde{\varphi}\!\left( \frac{s}{r}x, \frac{s}{r}x \right)
\end{eqnarray*}
for all $x \in {\mathcal X}$.
In a similar way we obtain the following inequality
\begin{eqnarray*}
\| h(x) - T(x) \| \leq \frac{1}{t}\, \varphi(0,x) + \frac{r}{t} \,
\widetilde{\varphi}\!\left( \frac{s}{r}x, \frac{s}{r}x \right)
\end{eqnarray*}
for all $x \in {\mathcal X}$.

If the mapping $\mu \mapsto f(\mu x)$ is continuous for each fixed
$x\in {\mathcal X}$, then the linearity of the mapping $T$ can be deduced by
the same reasoning as in the proof of the main theorem of
\cite{RAS1}.
\end{proof}

\begin{corollary} Suppose $f : {\mathcal X} \to {\mathcal Y}$ is a
mapping with $f(0) = 0$, and there exist constants
$\varepsilon, \delta \geq 0$ and $p \in [0, 1)$ such that
\begin{eqnarray*}
\left\| r f\!\left( \frac{s x + t y}{r} \right) - s f(x) - t f(y) \right\|
\leq \varepsilon + \delta ( \|x\|^p + \|y\|^p ),
\end{eqnarray*}
for all $x, y\in {\mathcal X}$.
Then there exists a unique additive mapping $T : {\mathcal X} \to {\mathcal Y}$
such that
\begin{eqnarray*}
\| f(x) - T(x) \| \leq \frac{3}{r}\varepsilon +
\frac{1}{r}\left[ \left( \frac{r}{s} \right)^p +
\left( \frac{r}{t} \right)^p \, \right]
\frac{2\delta \|x\|^p}{1-2^{p-1}}
\end{eqnarray*}
for all $x \in {\mathcal X}$.
\end{corollary}

\begin{proof} Define $\varphi(x, y) = \varepsilon + \delta (\|x\|^p + \|y\|^p)$,
and apply Theorem 2.1.
\end{proof}

\section{Asymptotic Behavior of Generalized Jensen Equation}

We start this section with investigating the stability of
generalized Jensen equation outside a ball.
The results in this section are generalizations of Theorem 3.3
and Corollary 3.4 of \cite{JUNG}.

\begin{theorem}
Let $d > 0$, $\varepsilon > 0$, and $f : {\mathcal X} \to {\mathcal Y}$
be a mapping with $f(0) = 0$ such that
\begin{eqnarray*}
\left\| r f\!\left( \frac{sx+ty}{r} \right) - s f(x) - t f(y) \right\|
\leq \varepsilon,
\end{eqnarray*}
for all $x, y \in {\mathcal X}$ with $\|x\| + \|y\| \geq d$.
Then there exists a unique additive mapping $T: {\mathcal X} \to {\mathcal Y}$
such that
\begin{eqnarray*}
\| f(x) - T(x) \| \leq \frac{15}{r}\, \varepsilon,
\end{eqnarray*}
for all $x \in {\mathcal X}$.
\end{theorem}

\begin{proof}
Let $x, y \in {\mathcal X}$ with $\|x\| + \|y\| < d$.
If $x = y = 0$, we choose a $z \in {\mathcal X}$ with $\| z \| = d$,
otherwise
\begin{eqnarray*}
z : = \left\{\begin{array}{ll}
                \left (1 + \frac{d}{\|x\|}\right) x & \mbox{if $\|x\| \geq \|y\|$} \medskip \\
               \left  (1 + \frac{d}{\|y\|}\right) y & \mbox{if $\|x\| \leq \|y\| . $}
             \end{array}
      \right.
\end{eqnarray*}
Then one can easily verify the following inequalities:
\begin{eqnarray*}
&&\left\| \left( 2 + \frac{t}{s} \right)z + \frac{t}{s}y \right\| +
  \left\| \frac{s}{t}x - \left( 1 + \frac{2s}{t} \right)z \right\| \geq d,  \\
&&\|x\| + \|z\| \geq d, \\
&&\left\| 2\left( 1 + \frac{t}{s} \right)z \right\| + \|y\| \geq d,  \\
&&\left\| 2\left( 1 + \frac{t}{s} \right)z \right\| +
  \left\| \frac{s}{t}x - \left( 1 + \frac{2s}{t} \right)z \right\| \geq d, \\
&&\left\| \left( 2 + \frac{t}{s} \right)z + \frac{t}{s}y \right\| + \|z\| \geq d.
\end{eqnarray*}
It follows that
\begin{eqnarray*}
\lefteqn{\left\| r f\!\left( \frac{sx+ty}{r} \right) - s f(x) - t f(y) \right\|}\\
& \leq & \left\| r f\!\left( \frac{sx+ty}{r} \right) -
                 s f\!\left( \left( 2 + \frac{t}{s} \right)z + \frac{t}{s}y \right) -
                 t f\!\left( \frac{s}{t}x - \left( 1 + \frac{2s}{t} \right)z \right)
         \right\| \\
&      & +\,\left\| r f\!\left( \frac{sx+tz}{r} \right) - s f(x) - t f(z) \right\| \\
&      & +\,\left\| r f\!\left( \frac{2(s+t)z+ty}{r} \right) -
                    s f\!\left( 2\left( 1 + \frac{t}{s} \right)z \right) - t f(y)
            \right\| \\
&      & +\,\left\| -r f\!\left( \frac{sx+tz}{r} \right) +
                    s f\!\left( 2\left( 1 + \frac{t}{s} \right)z \right) +
                    t f\!\left( \frac{s}{t}x - \left( 1 + \frac{2s}{t} \right)z \right)
            \right\| \\
&      & +\,\left\| -r f\!\left( \frac{2(s+t)z+ty}{r} \right) +
                    s f\!\left( \left( 2 + \frac{t}{s} \right)z + \frac{t}{s}y \right) +
                    t f(z)
            \right\| \\
& \leq & 5\varepsilon.
\end{eqnarray*}
Hence
\begin{eqnarray*}
\left\| r f\!\left( \frac{sx+ty}{r} \right) - s f(x) - t f(y) \right\|
\leq 5\varepsilon
\end{eqnarray*}
holds for all $x,y\in {\mathcal X}$.
Using Corollary 2.2 (with $\delta = 0$), we conclude the existence of a
unique additive mapping $T : {\mathcal X} \to {\mathcal Y}$ such that
\begin{eqnarray*}
\| f(x) - T(x) \| \leq \frac{15}{r}\varepsilon,
\end{eqnarray*}
for all $x \in {\mathcal X}$.
\end{proof}

Now we are ready to study the asymptotic behavior of generalized
Jensen equation.

\begin{corollary}
Let $f : {\mathcal X} \to {\mathcal Y}$ be a mapping with
$f(0)=0$. Then $f$ is additive if and only if
\begin{eqnarray}\label{asym}
\left \|r f \left (\frac{sx+ty}{r} \right ) - s f(x) - t f(y) \right \|\to 0 \; \; \; \; as \;
\; \; \; \|x\| + \|y\|\to \infty.
\end{eqnarray}
\end{corollary}
\begin{proof}
Let (\ref{asym}) be satisfied. Then there is a sequence
$\{\delta_n\}$ monotonically decreasing to zero such that
\begin{eqnarray}\label{deltan}
\left \|r f \left (\frac{sx+ty}{r}\right ) - s f(x) - t f(y)\right \|\leq \delta_n,
\end{eqnarray}
for all $x,y\in {\mathcal X}$ with $\|x\|+\|y\|\geq n$. Applying
(\ref{deltan}) and Theorem 3.1 we obtain a sequence $\{T_n\}$ of
unique additive mappings from ${\mathcal X}$ into ${\mathcal Y}$
such that
\begin{eqnarray*}
\|f(x) - T_n(x)\| \leq \frac{15}{r} \, \delta_n,
\end{eqnarray*}
for all $x\in {\mathcal X}$. The uniqueness of $T_n$ implies that
$T_n=T_{n+j}$ for all $j\in \N$. Hence by letting $n$ tend to
$\infty$ we infer that $f$ is additive. The reverse statement is
obvious.
\end{proof}

\section{Stability on the Punctured Space}

In this section we prove the stability of Jensen equation of
Pexider type on the deleted space ${\mathcal X}_0 : = {\mathcal
X} \setminus \{0\}$. In particular, in the case of $s = t = r$, we obtain
some results on stability of the Pexiderized Cauchy equation
restricted to the punctured space.

\begin{proposition}
Let $f, g, h : {\mathcal X} \to {\mathcal Y}$ be mappings with
$f(0) = g(0) = h(0) = 0$ for which there exists a function
$\varphi : {\mathcal X}_0\times {\mathcal X}_0 \to [0, \infty)$
satisfying
\begin{eqnarray*}
\widetilde{\varphi}(x, y) &: =& \frac{2}{3}\, \sum_{n=0}^{\infty} 3^{-n}
\bigg [ \varphi \left (\frac{3^{n+1}}{2}x,\frac{-3^n}{2}y \right )+ \frac{1}{2}
\varphi \left (\frac{3^{n+1}}{2}x,\frac{3^{n+1}}{2}y \right )\\
&&+ \frac{1}{2}\varphi \left (\frac{3^{n+1}}{2}x,\frac{-3^{n+1}}{2}y \right ) +
\frac{1}{2} \varphi \left (\frac{3^n}{2}x,\frac{3^n}{2}y \right ) + \frac{1}{2}
\varphi \left (\frac{3^n}{2}x,\frac{-3^n}{2}y \right )\bigg ]\\
&<& \infty
\end{eqnarray*}
and
\begin{eqnarray}\label{oJen}
\left \|r \, f \left (\frac{sx+ty}{r}\right ) - s \, g(x) - t \, h(y) \right \|\leq
\varphi \left (x,\frac{t}{s}y \right )
\end{eqnarray}
for all $x, y\in {\mathcal X}_0$. If $h$ is an odd
mapping, then there exists a unique additive mapping $A :
{\mathcal X}\to {\mathcal Y}$ such that
\begin{eqnarray*}
\|f(x) - A (x)\| \leq \frac{1}{r}
\widetilde{\varphi} \left (\frac{r}{s}x, \frac{r}{s}x \right ),
\end{eqnarray*}
\begin{eqnarray*}
\left \|g(x) - A(x) \right \| \leq \frac{1}{2s}\, \varphi(x,x) +
\frac{1}{2s}\, \varphi(x,-x) +\frac{1}{2s}\, \widetilde{\varphi}(2x,2x)
\end{eqnarray*}
and
\begin{eqnarray*}
\left \|h(x) - A(x) \right \| \leq
\frac{1}{2t}\, \varphi \left (\frac{t}{s}x,\frac{t}{s}x \right ) +
\frac{1}{2t}\, \varphi \left (\frac{t}{s}x,\frac{-t}{s}x \right )
+\frac{1}{2t}\, \widetilde{\varphi}\left (\frac{2t}{s}x,\frac{2t}{s}x \right )
\end{eqnarray*}
for all $x\in {\mathcal X}$.
\end{proposition}
\begin{proof}
Replacing $y$ by $\frac{s}{t}y$ in (\ref{oJen}) we get
\begin{eqnarray*}
\left \|r \, f \left (\frac{s}{r} \, (x+y) \right ) - s \, g(x)
- t \, h \left (\frac{s}{t}y \right )\right \|\leq
\varphi (x,y)
\end{eqnarray*}
for all $x, y\in {\mathcal X}$. Define the mappings $F, G, H$ by
$F(x) : = r \, f \left (\frac{s}{r}x \right ), G(x) : = s\, g(x)$ and $H(x) : = t \,
h\left (\frac{s}{t}y \right )$. Then
\begin{eqnarray}\label{olin}
\|F(x+y)) - G(x) - H(y)\|\leq \varphi(x,y)
\end{eqnarray}
for all $x, y\in {\mathcal X}_0$.

Replacing both $x$ and $y$ by $\frac{x}{2}$ in (\ref{olin}) we get
\begin{eqnarray}\label{1}
\left \|F(x) - G \left (\frac{x}{2} \right ) - H \left (\frac{x}{2} \right ) \right \| \leq
\varphi \left (\frac{x}{2},\frac{x}{2} \right )
\end{eqnarray}
for all $x\in {\mathcal X}_0$. Replacing $x$ by $\frac{x}{2}$ and
$y$ by $\frac{-x}{2}$, respectively, in (\ref{olin}) we obtain
\begin{eqnarray}\label{2}
\left \|G \left (\frac{x}{2} \right ) - H \left (\frac{x}{2}\right ) \right \| \leq
\varphi\left (\frac{x}{2}, \frac{-x}{2}\right )
\end{eqnarray}
for all $x\in {\mathcal X}_0$. It follows from (\ref{1}) and
(\ref{2}) that
\begin{eqnarray}\label{3}
\left \|F(x) - 2 \, H \left (\frac{x}{2}\right )\right \|
\leq \varphi\left (\frac{x}{2},\frac{x}{2}\right )+
\varphi\left (\frac{x}{2},\frac{-x}{2}\right )
\end{eqnarray}
and
\begin{eqnarray}\label{4}
\left \|F(x) - 2 \, G\left (\frac{x}{2}\right )\right \|
\leq \varphi\left (\frac{x}{2},\frac{x}{2}\right )+
\varphi\left (\frac{x}{2},\frac{-x}{2}\right )
\end{eqnarray}
for all $x\in {\mathcal X}_0$. Replacing $x$ by $\frac{3x}{2}$
and $y$ by $\frac{-x}{2}$, respectively, in (\ref{olin}) we have
\begin{eqnarray}\label{5}
\left \|F(x) - G\left (\frac{3x}{2}\right ) + H\left (\frac{x}{2} \right )\right \| \leq
\varphi \left (\frac{3x}{2},\frac{-x}{2}\right )
\end{eqnarray}
for all $x\in {\mathcal X}_0$. Using (\ref{3}), (\ref{4}) and
(\ref{5}) we infer that
\begin{eqnarray}\label{6}
\left \|F(x) - \frac{1}{3}F(3x)\right \|
&=& \frac{2}{3}\, \left \|F(x) -
\frac{1}{2}\, F(3x) +
\frac{1}{2}\, F(x) \right \|\nonumber\\
&\leq&   \frac{2}{3}\,
\left \|F(x) - G\left (\frac{3x}{2}\right ) + H\left (\frac{x}{2}\right )\right \|
+  \frac{2}{3}\, \left \|G \left (\frac{3x}{2}\right ) - \frac{1}{2}\, F(3x)\right \|\nonumber\\
&& +  \frac{2}{3}\,\left \|\frac{1}{2}\, F(x) - H\left (\frac{x}{2}\right )\right \|\nonumber\\
&\leq&  \frac{2}{3}\,\varphi\left (\frac{3x}{2},\frac{-x}{2}\right )
+\frac{1}{3}\, \varphi\left (\frac{3x}{2},\frac{3x}{2}\right )
+ \frac{1}{3}\, \varphi\left (\frac{3x}{2},\frac{-3x}{2}\right )\nonumber\\
&& + \frac{1}{3}\, \varphi\left (\frac{x}{2},\frac{x}{2}\right )
+ \frac{1}{3}\, \varphi\left (\frac{x}{2},\frac{-x}{2}\right )
: = \psi(x)
\end{eqnarray}
for all $x\in {\mathcal X}_0$. Using (\ref{6}) and the induction,
one can prove that
\begin{eqnarray}\label{8}
\left \|F(x) - \frac{1}{3^n}\, F\left (3^nx \right ) \right \| \leq
\sum_{k=0}^{n-1}3^{-k}\, \psi \left (3^kx \right )
\end{eqnarray}
for all $n$ and all $x\in {\mathcal X}_0$, and also
\begin{eqnarray*}
\left \|\frac{1}{3^n} \, F \left (3^nx \right )
- \frac{1}{3^m}\, F\left (3^mx \right ) \right \| \leq
\sum_{k=m}^{n-1}3^{-k}\, \psi \left (3^kx \right )
\end{eqnarray*}
for all $m<n$ and all $x \in {\mathcal X}_0$. Since ${\mathcal
Y}$ is complete we deduce that the sequence
$\left \{\frac{1}{3^n}F(3^nx)\right \}$ is convergent.
Therefore we can define
the mapping $T: {\mathcal X} \to {\mathcal Y}$ by
\begin{eqnarray}\label{9}
T(x) : =\lim_{n\to\infty}\, \frac{1}{3^n}F\left (3^nx\right )
= \lim_{n\to\infty}\frac{r}{3^n}\, f\left (\frac{3^n \, s}{r}x \right ).
\end{eqnarray}
It follows from (\ref{8}) and (\ref{9}) that
\begin{eqnarray}\label{fT}
\left \|r f\left (\frac{s}{r}x\right ) - T(x)\right \| = \|F(x) - T(x)\| \leq
\widetilde{\varphi}(x,x)
\end{eqnarray}
and so
\begin{eqnarray}\label{addeq}
\left \|f(x) - \frac{1}{r}\, T\left (\frac{r}{s} x \right )\right \| \leq \frac{1}{r}\,
\widetilde{\varphi}\left (\frac{r}{s}x, \frac{r}{s}x \right )
\end{eqnarray}
for all $x\in {\mathcal X}$. Note that $T(0) = 0$.

Using (\ref{olin}), (\ref{3}) and (\ref{4}) we have
\begin{eqnarray*}
\left\| 2 \, F\!\left( \frac{x+y}{2} \right) - F(x) - F(y) \right\|
& \leq & \varphi\!\left( \frac{x}{2}, \frac{x}{2} \right) +
         \varphi\!\left( \frac{x}{2}, \frac{-x}{2} \right) \\
&      & +\,2\varphi\!\left( \frac{x}{2}, \frac{y}{2} \right) +
         \varphi\!\left( \frac{y}{2}, \frac{y}{2} \right) +
         \varphi\!\left( \frac{y}{2}, \frac{-y}{2} \right)
\end{eqnarray*}
for all $x, y \in {\mathcal X}_0$, whence
\begin{eqnarray*}
\lefteqn{\left\| \frac{2}{3^n}F\!\left( \frac{3^n(x+y)}{2} \right) -
                 \frac{1}{3^n}F(3^n x) - \frac{1}{3^n}F(3^n y)
         \right\|} \qquad\quad\\
& \leq & \frac{1}{3^{n}}\varphi\!\left( \frac{3^nx}{2}, \frac{3^nx}{2} \right) +
         \frac{1}{3^{n}}\varphi\!\left( \frac{3^nx}{2}, \frac{-3^nx}{2} \right) +
         \frac{2}{3^{n}}\varphi\!\left( \frac{3^nx}{2}, \frac{3^ny}{2} \right) \\
&      & +\,\frac{1}{3^{n}}\varphi\!\left( \frac{3^ny}{2}, \frac{3^ny}{2} \right) +
         \frac{1}{3^{n}}\varphi\!\left( \frac{3^ny}{2}, \frac{-3^ny}{2} \right)
\end{eqnarray*}
for all $x\in {\mathcal X}_0$. Letting $n\to\infty$ and noting to
the fact that the right hand side tends to zero, we conclude that
$T$ satisfies the Jensen equation and so it is additive. By a
known strategy one can easily establish the uniqueness of $T$; cf.
\cite{MOS2}. It follows from (\ref{addeq}) that
$$
\left\| f(x) - \frac{1}{s}\, T(x) \right\| \leq
\frac{1}{r}\, \widetilde{\varphi}\!\left( \frac{r}{s}x, \frac{r}{s}x \right)
$$
for all $x \in {\mathcal X}$.

Using(\ref{4}) and (\ref{fT}), we obtain
\begin{eqnarray*}
\|s g(x) - T(x)\| = \|G(x) - T(x)\|&\leq& \left \|G(x) -
\frac{1}{2}F(2x)\right \|
+\left \|\frac{1}{2}F(2x) - \frac{1}{2}T(2x)\right \|\\
&\leq& \frac{1}{2}\, \varphi(x,x) + \frac{1}{2}\, \varphi(x,-x)
+\frac{1}{2}\, \widetilde{\varphi}(2x,2x),
\end{eqnarray*}
and so
\begin{eqnarray*}
\left \|g(x) - \frac{1}{s} T(x)\right \| \leq \frac{1}{2s}\, \varphi(x,x) +
\frac{1}{2s}\, \varphi(x,-x) +\frac{1}{2s}\, \widetilde{\varphi}(2x,2x)
\end{eqnarray*}
for all $x\in {\mathcal X}$. Similarly by applying (\ref{fT}) and
(\ref{3}) we have
\begin{eqnarray*}
\left \|t \, h\left (\frac{s}{t}x \right ) - T(x)\right \|
= \|H(x) - T(x)\|\leq \frac{1}{2}\, \varphi(x,x) + \frac{1}{2}\, \varphi(x,-x)
+\frac{1}{2}\, \widetilde{\varphi}(2x,2x),
\end{eqnarray*}
and so
\begin{eqnarray*}
\left \|h(x) - \frac{1}{s} T(x) \right \| \leq
\frac{1}{2t}\, \varphi\left (\frac{t}{s}x,\frac{t}{s}x \right ) +
\frac{1}{2t}\, \varphi\left (\frac{t}{s}x,\frac{-t}{s}x\right )
+\frac{1}{2t}\, \widetilde{\varphi}\left (\frac{2t}{s}x,\frac{2t}{s}x\right )
\end{eqnarray*}
for all $x\in {\mathcal X}$.

With $A(x) := \frac{1}{s}T(x)$, we conclude that our assertions are true.
\end{proof}

\begin{proposition}
Let $f, g, h : {\mathcal X} \to {\mathcal Y}$ be mappings with
$f(0) = g(0) = h(0) = 0$ for which there exists a function
$\varphi : {\mathcal X}_0\times {\mathcal X}_0 \to [0, \infty)$
such that
\begin{eqnarray}\label{eJen}
\left \|r\, f\left (\frac{sx+ty}{r}\right ) - s\, g(x) - t \, h(y) \right \|\leq
\varphi\left (x,\frac{t}{s}y\right )
\end{eqnarray}
for all $x, y\in {\mathcal X}_0$. If $h$ is an even
mapping, then
\begin{eqnarray*}
\|f(x)\| \leq \frac{1}{r}\, \varphi\left (\frac{r x}{2s},\frac{r x}{2s}\right )
+\frac{1}{r}\, \varphi\left (\frac{r x}{2s},\frac{-r x}{2s}\right ),
\end{eqnarray*}
and
\begin{eqnarray*}
\left \|g(x) - \frac{t}{s}\,  h \left (\frac{s}{t} x\right ) \right \| \leq \frac{1}{s} \,
\varphi(x,-x)
\end{eqnarray*}
for all $x\in {\mathcal X}$.
\end{proposition}
\begin{proof}
Replacing $y$ by $\frac{s}{t}y$ in (\ref{eJen}) we get
\begin{eqnarray*}
\left \|r \, f\left (\frac{s}{r} (x+y) \right ) - s g(x)
- t \, h\left (\frac{s}{t}y\right )\right \|\leq
\varphi (x,y)
\end{eqnarray*}
for all $x, y\in {\mathcal X}$. Define the mappings $F, G, H$ by
$F(x) : = r \, f\left (\frac{s}{r}x\right ), G(x) : = s \, g(x)$ and $H(x) : = t\,
h\left (\frac{s}{t}y\right )$. Then
\begin{eqnarray}\label{elin}
\|F(x+y)) - G(x) - H(y)\|\leq \varphi(x,y)
\end{eqnarray}
for all $x, y\in {\mathcal X}_0$.

Replacing both $x$ and $y$ by $\frac{x}{2}$ in (\ref{elin}) we get
\begin{eqnarray}\label{10}
\left \|F(x) - G\left (\frac{x}{2}\right ) - H\left (\frac{x}{2}\right )\right \| \leq
\varphi\left (\frac{x}{2},\frac{x}{2}\right )
\end{eqnarray}
for all $x\in {\mathcal X}_0$. Replacing $x$ by $\frac{x}{2}$ and
$y$ by $\frac{-x}{2}$, respectively, in (\ref{elin}) we obtain
\begin{eqnarray}\label{11}
\left \|G\left (\frac{x}{2}\right ) + H\left (\frac{x}{2}\right )\right \| \leq
\varphi\left (\frac{x}{2},\frac{-x}{2}\right )
\end{eqnarray}
for all $x\in {\mathcal X}_0$. It follows from (\ref{10}) and
(\ref{11}) that
\begin{eqnarray*}
\left \|r \, f\left (\frac{s}{r}x\right )\right \|
= \|F(x)\|&\leq& \left \|F(x) - G\left (\frac{x}{2}\right ) -
H\left (\frac{x}{2}\right )\right \| +
\left \|G\left (\frac{x}{2}\right ) + H\left (\frac{x}{2}\right )\right \|\\
&\leq& \varphi\left (\frac{x}{2},\frac{x}{2}\right ) +
\varphi\left (\frac{x}{2},\frac{-x}{2}\right ),
\end{eqnarray*}
and so
\begin{eqnarray*}
\|f(x)\| \leq \frac{1}{r}\, \varphi\left (\frac{r x}{2s},\frac{r x}{2s}\right ) +
\frac{1}{r}\, \varphi\left (\frac{r x}{2s},\frac{-r x}{2s}\right ).
\end{eqnarray*}
In addition
\begin{eqnarray*}
\left \|g(x) - \frac{t}{s} \, h\left (\frac{s}{t}x \right ) \right \|
= \frac{1}{s}\, \|G(x) + H(x)\|
\leq \frac{1}{s}\, \varphi(x,-x)
\end{eqnarray*}
for all $x\in {\mathcal X}$.
\end{proof}

\begin{theorem}
Let $\varepsilon > 0$ and let $f : {\mathcal X} \to {\mathcal Y}$ be a
mapping with $f(0) = 0$ satisfying
\begin{eqnarray}\label{eps}
\left\| r \, f\!\left( \frac{sx+ty}{r} \right) - s \, f(x) - t \, f(y) \right\|
\leq \varepsilon
\end{eqnarray}
for all $x, y \in {\mathcal X}_0$.
Then there exists a unique additive mapping $A : {\mathcal X}\to {\mathcal Y}$
such that
\begin{eqnarray*}
\| f(x) - A(x) \|
\leq \min\!\left\{ \frac{3\varepsilon}{r}, \frac{5\varepsilon}{2s},
                   \frac{5\varepsilon}{2t} \right\} + \frac{2\varepsilon}{r}
\end{eqnarray*}
for all $x \in {\mathcal X}$.
\end{theorem}
\begin{proof}
Replacing $x, y$ by $-x, - \frac{sy}{t}$ in (\ref{eps}),
respectively, we get
\begin{eqnarray}\label{mineps}
\left \|r \, f\left (\frac{s}{r} (-x-y)\right )
- s \, f(-x) - t \, f\left (-\frac{s}{t}y\right )\right \|\leq
\epsilon
\end{eqnarray}
for all $x, y\in {\mathcal X}_0$. Using the odd and even parts
$f^o, f^e$ of $f$ and (\ref{eps}) and (\ref{mineps}) we obtain
\begin{eqnarray*}
\left \|r \, f^o\left (\frac{sx+ty}{r}\right ) - s \, f^o(x)
- t \, f^o(y)\right \|\leq \epsilon
\end{eqnarray*}
and
\begin{eqnarray*}
\left \|r \, f^e\left (\frac{sx+ty}{r}\right ) - s\, f^e(x)
- t \,f^e(y)\right \|\leq \epsilon
\end{eqnarray*}
for all $x\in {\mathcal X}_0$. Then Propositions 4.1 and 4.2 give
us a unique additive mapping $A$ such that
\begin{eqnarray*}
\left \|f^o(x) - A(x)\right \| \leq
\min\!\left\{ \frac{3\varepsilon}{r}, \frac{5\varepsilon}{2s},
                   \frac{5\varepsilon}{2t} \right\}
\end{eqnarray*}
and
\begin{eqnarray*}
\|f^e (x)\| \leq \frac{2\epsilon}{r}
\end{eqnarray*}
for all $x\in {\mathcal X}$. Hence
\begin{eqnarray*}
\|f(x) - A(x)\| \leq \|f^0 - A(x)\| + \|f^e\|\,
\min\!\left\{ \frac{3\varepsilon}{r}, \frac{5\varepsilon}{2s},
                   \frac{5\varepsilon}{2t} \right\}
                    + \frac{2\epsilon}{r}
\end{eqnarray*}
for all $x\in {\mathcal X}$.
\end{proof}

\section{Stability on Orthogonality Spaces}

Let us recall the orthogonality in the sense of R\" atz; cf. \cite{RAT}.
Suppose ${\mathcal X}$ is a real vector space with $\dim {\mathcal X}\geq 2$
and $\perp$ is a binary relation on ${\mathcal X}$ with the following properties:
\begin{itemize}
\item[{\rm (O1)}] $x\perp 0, 0\perp x$ for all $x\in {\mathcal X}$;
\item[{\rm (O2)}] if $x,y\in {\mathcal X}-\{0\}, x\perp y$, then $x,y$ are
                  linearly independent;
\item[{\rm (O3)}] if $x,y\in {\mathcal X}, x\perp y$, then
                  $\alpha x \perp \beta y$ for all $\alpha, \beta \in \R$;
\item[{\rm (O4)}] if $P$ is a $2$-dimensional subspace of ${\mathcal X}, x \in P$
                  and $\lambda \in \R_+$, then there exists $y_0 \in P$ such that
                  $x \perp y_0$ and $x + y_0 \perp \lambda x - y_0$.
\end{itemize}
The pair $({\mathcal X},\perp)$ is called an orthogonality space.
By an orthogonality normed space we mean an orthogonality space
having a normed structure.

Some interesting examples are as follows.
\begin{itemize}
\item[{\rm (i)}] The trivial orthogonality on a vector space ${\mathcal X}$
                 defined by (O1), and for non-zero elements
                 $x, y \in {\mathcal X}$, $x \perp y$ if and only if
                 $x, y$ are linearly independent.
\item[{\rm (ii)}] The ordinary orthogonality on an inner product space
                  $({\mathcal X}, \langle\,.,.\,\rangle)$ given by
                  $x \perp y$ if and only if $\langle x, y \rangle = 0$.
\item[{\rm (iii)}] The Birkhoff--James orthogonality on a normed space
                   $({\mathcal X}, \|.\|)$ defined by $x \perp y$ if and only if
                   $\| x + \lambda y \| \geq \|x\|$ for all $\lambda \in \R$.
\end{itemize}

Let ${\mathcal X}$ be a vector space (an orthogonality space) and
$({\mathcal Y},+)$ be an abelian group. Then a mapping $f:
{\mathcal X} \to {\mathcal Y}$ is called

\begin{itemize}
\item[{\rm (i)}] {\it orthogonally additive} if it satisfies the additive
                 functional equation for all $x, y \in {\mathcal X}$ with
                 $x \perp y$;
\item[{\rm (ii)}] {\it orthogonally quadratic} if it satisfies the quadratic
                  functional equation $f(x+y) + f(x-y) = 2f(x) + 2f(y)$
                  for all $x, y \in {\mathcal X}$ with $x\perp y$;
\item[{\rm (iii)}] {\it orthogonally generalized Jensen} if it satisfies
                   the generalized Jensen functional equation for all
                   $x, y \in {\mathcal X}$ with $x \perp y$.
\end{itemize}

In \cite{MOS2, MOS3} the orthogonal stability of various
functional equations were established. In particular, the
following theorem was proved; cf. Theorem 1 of \cite{MOS2}.

\begin{theorem}
Suppose $({\mathcal X},\perp)$ is an orthogonality space and
$({\mathcal Y}, \|.\|)$ is a real Banach space. Let $f, g,
h:{\mathcal X}\to {\mathcal Y}$ be mappings fulfilling
\begin{eqnarray*}
\|f(x+y)-g(x)-h(y)\|\leq\epsilon
\end{eqnarray*}
for some $\epsilon$ and for all $x,y\in {\mathcal X}$ with
$x\perp y$. Then there exists exactly a quadratic mapping
$Q:{\mathcal X}\to {\mathcal Y}$ and an additive mapping
$T:{\mathcal X}\to {\mathcal Y}$ such that
\begin{eqnarray*}
\|f(x)-f(0)-Q(x)-T(x)\|\leq \frac{68}{3}\epsilon , \\
\|g(x)-g(0)-Q(x)-T(x)\|\leq \frac{80}{3}\epsilon , \\
\|h(x)-h(0)-Q(x)-T(x)\|\leq \frac{80}{3}\epsilon
\end{eqnarray*}
for all $x\in {\mathcal X}$.
\end{theorem}
Now we are ready to provide another result on orthogonal
stability of functional equations.

\begin{theorem}
Suppose that $f, g, h : {\mathcal X}\to {\mathcal Y}$ is a
mapping with $f(0) = g(0) = h(0) = 0$ satisfying
\begin{eqnarray*}
\left \|r \, f \left (\frac{sx+ty}{r}\right ) - s \, g(x)
- t \, h(y)\right \|\leq\epsilon
\end{eqnarray*}
for all $x, y\in {\mathcal X}$ with $x\perp y$. Then there exist
a unique additive mapping $T: {\mathcal X}\to {\mathcal Y}$ and a
unique quadratic mapping $Q: {\mathcal X}\to {\mathcal Y}$ such
that
\begin{eqnarray*}
\|f(x) - T(x) - Q(x)\| \leq 68\epsilon ,
\end{eqnarray*}
\begin{eqnarray*}
\|g(x) - T(x) - Q(x)\| \leq 80\epsilon ,
\end{eqnarray*}
\begin{eqnarray*}
\|h(x) - T(x) - Q(x)\| \leq 80\epsilon
\end{eqnarray*}
for all $x\in {\mathcal X}$.
\end{theorem}
\begin{proof}
Using the same argument as in the proof of Theorem 2.1 we conclude
that $\|f(x+y) - g(x) - h(x)\|\leq 3\epsilon$. Then the result is
followed from Theorem 5.1.
\end{proof}

\section{Orthogonally generalized Jensen mappings on Balls in Inner Product Spaces}

Sikorska showed that if $f$ is an orthogonally additive mapping
on an open ball ${\bf B}$ of a real inner product space
${\mathcal X}$ into a real sequentially complete linear
topological space ${\mathcal Y}$ then there exist additive
mappings $T: {\mathcal X}\to {\mathcal Y}$ and $b:\R_+\to
{\mathcal Y}$ such that $f(x) = T(x) + b(\|x\|^2)$ for all $x\in
{\bf B}$; cf. Corollary 1 of \cite{SIK}. By an orthogonally
generalized Jensen mapping we mean a mapping $f : {\mathcal X}
\to {\mathcal Y}$ such that
$\left \|r \, f\left (\frac{sx+ty}{r}\right ) - s \, f(x) - t\,
f(y)\right \|\leq\epsilon$ holds for some $\epsilon>0$ and for all $x,
y\in {\mathcal X}$ with $x\perp y$. We can extend Sikorska's
result to the orthogonally generalized Jensen mappings as follows.
\begin{theorem}
If ${\bf B}$ is an open ball of a real inner product space
${\mathcal X}$, ${\mathcal Y}$ is a real sequentially complete
linear topological space, and $f:{\bf B}\to {\mathcal Y}$ is
orthogonally generalized Jensen with $f(0) = 0$, then there exist
additive mappings $T: {\mathcal X}\to {\mathcal Y}$ and $b:\R_+\to
{\mathcal Y}$ such that $f(x) = T(x) + b\left (\|x\|^2\right )$ for all $x\in
{\bf B}$.
\end{theorem}
\begin{proof} With $y=0$, we have
\begin{eqnarray*}
r \, f\left (\frac{sx+t0}{r}\right ) = s \, f(x) + t \, f(0) = s f(x),
\end{eqnarray*}
and therefore $f\left (\frac{r}{s} x\right ) = \frac{r}{s} f(x)$. Similarly
$f(\frac{r}{t} x) = \frac{r}{t} f(x)$. Hence
\begin{eqnarray*}
r \, f(x+y)=s\, f\left (\frac{r}{s} x\right ) + t \, f\left (\frac{r}{t} y\right )
= r \, f(x) + r\, f(y)
\end{eqnarray*}
and so $f(x + y) = f(x) + f(y)$ for all $x,y\in {\bf B}$ with $x \perp y$. This
implies that $f$ is orthogonally additive on ${\bf B}$. Now the
assertion follows from the Sikorska's result.
\end{proof}
If we are restricted to the punctured ball ${\bf B}\setminus\{0\}$ then by
following the same strategy as in Lemma 1 of \cite{SIK} we obtain
the following result in the case that $s=t>\frac{1}{\sqrt{2}} \, r$.
The situation in the general case left as a question.

\begin{theorem}
If ${\bf B}$ is an open ball of a real inner product space
${\mathcal X}$ of dimension greater than $1$,
${\mathcal Y}$ is a real sequentially complete
linear topological space, and $f: {\bf B}\setminus\{0\} \to {\mathcal Y}$
is orthogonally generalized Jensen map with parameters
$s=t>\frac{1}{\sqrt{2}} \, r$, then there exist additive mappings
$T: {\mathcal X}\to {\mathcal Y}$ and $b:\R_+\to {\mathcal Y}$
such that $f(x) = T(x) + b\left (\|x\|^2\right )$ for all $x\in {\bf B}\setminus \{0\}$.
\end{theorem}
\begin{proof}
First note that if $f$ is a generalized Jensen map with parameters
$t=s \geq r $, then $f(\lambda(x+y))=\lambda f(x) + \lambda f(y)$ for
some $\lambda \geq 1$ and all $x, y\in {\bf B}\setminus \{0\}$ such that $x \perp y$.

\medskip
\noindent
\underline{\rm Step (I)- the case that f is odd:}
Let $x \in {\bf B} \setminus \{0\}$. There exists $y_0 \in {\bf B} \setminus
\{0\}$ such that $x \perp y_0$, $x + y_0 \perp x - y_0$. Since
$\frac{x + y_0}{2 \lambda}, \frac{x - y_0}{2 \lambda}, \frac{x}{2
\lambda^2}, \frac{y_0}{2 \lambda^2} \in {\bf B}$, we have
\begin{eqnarray*}
f(x)&=& f(x)- \lambda\, f\left ( \frac{x+y_0}{2\, \lambda}\, \right ) -
\lambda \, f\left ( \frac{x-y_0}{2\, \lambda}\, \right )\\
&&+ \, \lambda \, f\left ( \frac{x+y_0}{2\, \lambda}\, \right ) - \lambda^2\,
f\left ( \frac{x}{2\, \lambda^2}\, \right ) - \lambda^2 \, f\left (
\frac{y_0}{2\, \lambda^2}\, \right )\\
&&+ \, \lambda \, f\left ( \frac{x-y_0}{2\, \lambda}\, \right ) - \lambda^2
f\left ( \frac{x}{2\, \lambda^2}\, \right ) - \lambda^2\,  f\left (
\frac{-y_0}{2\, \lambda^2}\, \right )\\
&&+\, 2\, \lambda^2 \, f\left (\frac{x}{2\, \lambda^2}\, \right )
= 2\, \lambda^2 \, f\left ( \frac{x}{2\, \lambda^2}\, \right ).
\end{eqnarray*}
It follows that
\begin{eqnarray}\label{lambda}
f(x)  =  2\lambda^2 f\!\left( (2\lambda^2)^{-1} \, x \right)
\end{eqnarray}
for all $x\in {\bf B} \setminus \{0\}$.
Since
\begin{eqnarray*}
\left (2\lambda^2\right )^{n+m}\, f\left ((2\lambda^2)^{-n-m}\, x\right )
= \left (2\lambda^2\right )^n \, f\left ((2\lambda^2)^{-n}\, x\right )
\end{eqnarray*}
for all $x\in {\bf B}\setminus \{0\}$, and all $m,n\in \N$ with
$(2\lambda^2)^{-n}x\in {\bf B}\setminus \{0\}$, we can well define the mapping
$T: {\mathcal X}\to {\mathcal Y}$ by $T(x) : = (2\lambda^2)^n\,
f((2\lambda^2)^{-n}\, x)$ where $n$ is an integer such that
$(2\lambda^2)^{-n}x\in {\bf B}\setminus \{0\}$.
Clearly, $T$ is an extension of $f$ to ${\mathcal X}$.
The mapping $T$ is an odd
orthogonally additive mapping. To see this, let $x,y\in {\mathcal X}$
with $x \perp y$.
Then there exists a positive integer $n$ such that $(2\lambda^2)^{-n}x, \,
(2\lambda^2)^{-n}y, \, (2\lambda^2)^{-n}(x+y), \,
(2\lambda^2)^{-n}(x-y)\in {\bf B}\setminus \{0\}$, and by using
 (\ref{lambda}) we obtain
\begin{eqnarray*}
T(x) + T(y) & = & (2\lambda^2)^{n} f((2\lambda^2)^{-n}x) +
                  (2\lambda^2)^n f((2\lambda^2)^{-n}y) \\
            & = & (2\lambda^2)^{n+1} f((2\lambda^2)^{-n-1}x) +
                  (2\lambda^2)^{n+1} f((2\lambda^2)^{-n-1}y)\\
            & = & (2\lambda^2)^{n+1} f((2\lambda^2)^{-n-1}(x +
            y))\\
            & = & (2\lambda^2)^{n} f((2\lambda^2)^{-n}(x + y))\\
            & = & T(x+y).
\end{eqnarray*}
Hence $T(x) + T(y) = T(x + y)$ for all $x, y \in {\mathcal X}$. By
Corollary 7 of \cite{RAT}, $T$ is additive and $f = T|_{{\bf B}\setminus\{0\}}$.

\medskip
\noindent
\underline{\rm Step (II)- the case that f is even:}
Using the same notation and the same reasoning as step (I), one
can show that $f(x)=f(y_0)$ and the mapping $Q: {\mathcal X}\to
{\mathcal Y}$ defined by $Q(x) : = (4\lambda^2)^n
f((2\lambda^2)^{-n}x)$ is even orthogonally additive. By
Corollaries 7 and 10 of \cite{RAT}, $Q$ is quadratic and there
exists an additive mapping $b:\R_+\to {\mathcal Y}$ such that
$Q(x) = b(\|x\|^2)$ for all $x\in {\bf B}$. In addition,
$f=Q|_{{\bf B}\setminus\{0\}}$

\medskip
\noindent
\underline{\rm Step (III)- the general case:}
If $f$ is an arbitrary mapping, then $f$ can be expressed as $f = f^o + f^e$
where $f^o$ and $f^e$ are the odd and even parts of $f$, respectively.
Now the result can be deduced from Steps (I) and (II).
\end{proof}

\end{document}